%multider.tex:
%There are EXACTLY 1493804444499093354916284290188948031229880469556 Ways to Derange a Standard  Deck of Cards (ignoring suits)
%[and many other such useful facts]
%%a Plain TeX file by Shalosh B. Ekhad, Christoph Koutschan and Doron Zeilberger (x pages)

%begin macros

\baselineskip=14pt
\parskip=10pt

\font\eightrm=cmr8 

\magnification=\magstephalf

\def\1{{\overline{1}}}
\def\2{{\overline{2}}}
\parindent=0pt
\overfullrule=0in

\def\frac#1#2{{#1 \over #2}}
%\headline={\rm  \ifodd\pageno  \RightHead  \else  \LeftHead  \fi}
%\def\RightHead{\centerline{
%Title
%}}
%\def\LeftHead{ \centerline{Doron Zeilberger}}
%end macros
\centerline
{\bf  
There are EXACTLY 1493804444499093354916284290188948031229880469556 
}
\centerline
{\bf Ways to Derange a Standard  Deck of Cards (ignoring suits) [and many other such useful facts]}
\bigskip
\centerline
{\it Shalosh B. EKHAD, Christoph KOUTSCHAN, and Doron ZEILBERGER}

\bigskip

\qquad \qquad {\it In fond memory of Joe Gillis (3 Aug. 1911- 18 Nov. 1993), who taught us that Special Functions Count }

The famous {\it probl\`eme des rencontres}, due to Montmort, can be formulated (see [R], p. 58) in terms of
counting the number of ways of rearranging a deck of $n$ {\bf distinct} cards, such that {\bf no} card landed in the previous location. 
Calling this number $D_n$, there are quite a few 
{\bf good answers}. Perhaps the best is as the {\bf solution} of the (inhomogeneous)  first-order linear recurrence
$$
D_{n+1} - (n+1) D_{n} = (-1)^{n+1} \quad, \quad D_1=0 \quad ,
$$
that enables you, {\it very fast}, to find the first, say, $10000$ terms.
It follows, in particular that the number of ways of completely {\it deranging} a standard deck of cards, where all the $52$ cards are
considered {\it distinct} is
$$
29672484407795138298279444403649511427278111361911893663894333196201 \quad .
$$

This raises the natural question, of counting {\it derangements} of a standard deck of cards with $13$ {\bf different} {\it denominations}
$\{1, \dots , 13\}$ (where $1=$Ace, $11=$Jack,  $12=$Queen,  $13=$King), but {\bf ignoring suits}. To our great surprise
this number was {\bf not} (Jan. 22, 2021) in Neil Sloane's monumental [OEIS]. The {\bf exact} value is the one in the title, namely 
$$
1493804444499093354916284290188948031229880469556 \quad .
$$

This is a special case of the problem of finding the number of {\bf derangements} of a {\bf multiset}.

You have $a_1$ copies of $1$, $a_2$ copies of $2$, $\dots$, $a_n$ copies of $n$, in other words you have the
{\bf multiset} $1^{a_1} \dots n^{a_n}$, and you are interested in the number, let's call it $M(a_1, \dots, a_n)$,
of {\bf derangements}. In other words, out of the total number of arrangements, given by the multinomial coefficient 
$(a_1+ \dots +a_n)!/(a_1! \cdots a_n!)$,
how many are there where each location is {\bf different} than the original.
Here is why such  questions are {\bf useful}.\footnote
{${}^1$}  
{ \eightrm
This answers a question of Persi Diaconis , who asked `Why is it useful?', at the end of DZ's talk about this problem,
at the special session of the Joint Mathematics Meeting in memory of Richard Askey, Jan. 2021. We thank
him for raising this important question, that also inspired the word `useful' in the title.}

Suppose that in a certain army there are $n$ different ranks, and each officer must wear a jacket 
indicating his or her rank. In a certain party, there were $a_1$ officers of rank $1$, $a_2$ officers of rank $2$,
$\dots$, $a_n$ officers of rank $n$. Before the party started, they hung their jackets  in the coat-room.
When the party was over, they were all very drunk, so they each grabbed a jacket (uniformly) at random.

{\it What is the probability that every one shows the wrong rank? }

Knowing the {\bf exact} (rather than the {\it approximate}) probability
is (potentially) {\bf very useful} in deciding what odds to bet on such an event. Recall that probability theory started out
with such {\bf important} (and useful!) questions about gambling!

For example, $M(2,2,2)=10$, since there are ten of them. Here they are:
$$
\left ( \matrix{ 1 & 1 & 2 & 2 & 3 & 3 \cr
                 2 & 2 & 3 & 3 & 1 & 1}
  \right ) \quad , 
\left ( \matrix{ 1 & 1 & 2 & 2 & 3 & 3 \cr
                 2 & 3 & 1 & 3 & 1 & 2}
  \right ) ,
\quad
\left ( \matrix{ 1 & 1 & 2 & 2 & 3 & 3 \cr
                 2 & 3 & 1 & 3 & 2 & 1}
  \right ) ,
\quad
\left ( \matrix{ 1 & 1 & 2 & 2 & 3 & 3 \cr
                 2 & 3 & 3 & 1 & 1 & 2}
  \right ) ,
$$
$$
\quad
\left ( \matrix{ 1 & 1 & 2 & 2 & 3 & 3 \cr
                 2 & 3 & 3 & 1 & 2 & 1}
  \right ) ,
\quad
\left ( \matrix{ 1 & 1 & 2 & 2 & 3 & 3 \cr
                 3 & 2 & 1 & 3 & 1 & 2}
  \right ) ,  \quad
\left ( \matrix{ 1 &  1 & 2 & 2 & 3 & 3 \cr
                 3  & 2 & 1 & 3 & 2 & 1}
  \right ) ,  \quad
\left ( \matrix{ 1 & 1 & 2 & 2 & 3 & 3 \cr
                 3 & 2 & 3 & 1 & 1 & 2}
  \right ) ,  \quad
$$
$$
\quad
\left ( \matrix{ 1 & 1 & 2 & 2 & 3 & 3 \cr
                 3 & 2 & 3 & 1 & 2 & 1}
  \right ) ,  \quad
\left ( \matrix{ 1 & 1 & 2 & 2 & 3 & 3 \cr
                 3 & 3 & 1 & 1 & 2 & 2}
  \right )   \quad .
$$

The million dollar question is: How to {\it efficiently} compute these numbers $M(a_1, a_2, \dots, a_n)$?
(for example, the number of the title is $M(4,4,\dots,4)$ where the $4$ is repeated $13$ times.)

For a fixed $n$, the great Percy Alexander MacMahon gave a nice {\it generating function} that he deemed important enough to
include in the entry {\it Combinatory Analysis} that he contributed, in 1910, to the eleventh edition of 
{\it Encyclopedia Britannica} ([M], p. 755, col. 1)
$$
\sum_{0 \leq a_1,\dots, a_n < \infty}  M(a_1, \dots, a_n)\, x_1^{a_1} \cdots x_n^{a_n} \, = \,
\frac{1}
{1-e_2-2e_3- \dots - (n-1) e_n} \quad,
$$
where $e_i=e_i(x_1, \dots, x_n)$ is the {\it elementary symmetric function} of degree $i$.

But this is not very efficient if one wants to find, say, $M(5,\dots, 5)$, where $5$ is repeated $100$ times.
One would have to find the Taylor series of a rational function with $100$ variables, and extract the coefficient
of $x_1^5 \cdots  x_{100}^5$. Even for the number in the title, one would have to extract the coefficient of $x_1^4 \cdots x_{13}^4$ in
a rational function with $13$ variables.

A much more efficient formula was given in the seminal paper [EG], by Shimon Even and Joe Gillis. It is the following
amazing formula, made more widely known by Richard Askey in his classic monograph ([A], p. 43).

{\bf Theorem} (Even-Gillis [EG]) The number of ways of deranging the multiset $1^{a_1} \dots n^{a_n}$ is
$$
(-1)^{a_1 + \dots +a_n} \, \int_0^{\infty} e^{-x} L_{a_1}(x) \, \cdots L_{a_n}(x) \, dx \quad  ,
$$
where $L_a(x)$ is the (simple) {\bf Laguerre polynomial} of degree $a$, so useful in classical potential theory, {\it Quantum Mechanics},
and elsewhere, and thanks to Gillis and Even, even in  combinatorics! 
$$
L_a(x) \, := \, \sum_{\alpha=0}^a \, (-1)^\alpha {{a} \choose {\alpha}} \frac{x^{\alpha}}{\alpha!} \quad.
$$

In particular the number of ways of deranging a standard deck of cards (given in the title) is
$$
\int_0^{\infty} e^{-x} L_4(x)^{13} \, dx \quad ,
$$
and it took Maple {\tt 0.123} seconds to find it! Using MacMahon's formula would take much longer!

From now on, we are interested in fast efficient computations of the numbers [$k$ repeated $n$ times below].
$$
F[n](k)=M(k,k,\dots , k)=(-1)^{kn} \int_0^{\infty} \, e^{-x} L_k(x)^n \, dx \quad .
$$

Suppose that we have $1000$ different denominations (rather than $13$) and still $4$ suits, then we would need the number $F[1000](4)$.
Now Maple has to integrate
$e^{-x}$ times a polynomial of degree $4000$ and it takes {\bf much longer!}, in fact, SBE took 295 seconds to compute it.

It so happens that thanks to the  {\it Holonomic Systems Approach} to {\bf Special Functions} [Z], initiated in 1990 by one of us (DZ),
and later extended and efficiently implemented [K] by another one of us (CK),
the Even-Gillis formula can be used to compute very fast $F[n](k)$ for small (and not so small) $k$, and very large $n$,
and also for small $n$ (we went as far as $n=9$) and very large $k$, for which doing it directly would take a very long time.

Recall that a function $f(n,x)$ is holonomic in $(n,x)$ if (roughly) it satisfies a {\bf linear differential equation} 
(with respect to $x$) with coefficients
that are polynomials in $(n,x)$ and a {\bf linear difference equation} (also called {\it recurrence equation}),
with respect to $n$, also with coefficients
that are polynomials in $(n,x)$. This is obviously true for the Laguerre polynomials $L_n(x)$.

It is easy to see 
that the integrand  $e^{-x} L_k(x)^n$ is holonomic in $(n,x)$ for each {\bf specific} (numeric) $k$.
In fact this is true for any function of $(n,x)$ that has the form $e^{-x} P(x)^n$, for an {\it arbitrary} polynomial $P(x)$.

It is also easy, since, as observed  in [Z], the class of holonomic functions is an {\bf algebra},
that the integrand  $e^{-x} L_k(x)^n$ is holonomic in $(k,x)$ for each {\bf specific} (numeric) $n$.

For the former case things are simpler, since the integrand is more than `just holonomic', is it {\bf hyper-exponential},
i.e. the two relevant equations (differential and difference) are {\bf first order}. For this special case
the {\bf Almkvist-Zeilberger algorithm} [AlZ]  implemented in the Maple package

{\tt https://sites.math.rutgers.edu/\~{}zeilberg/tokhniot/EKHAD.txt} \quad 

(procedure AZdI), and also included in this article's  Maple package 

{\tt https://sites.math.rutgers.edu/\~{}zeilberg/tokhniot/MultiDer.txt} \quad  ,

finds such a recurrence very fast. 

The output file
{\tt https://sites.math.rutgers.edu/\~{}zeilberg/tokhniot/oMultiDer1.txt} \quad 
contains such recurrences for $1 \leq k \leq 19$.

The [OEIS] (viewed Jan. 21, 2021) only has the sequences for $1 \leq k \leq 5$.

The case $k=1$ is the classical derangement sequence  {\tt A000166}. The case $k=2$ is {\tt A000459}, where
a recurrence is given. The case $k=3$  is {\tt A059073}, that only has the entries up to $n=12$, and there is no recurrence.
The case $k=4$  is {\tt A059074}, where also only the entries up to $n=12$ are given, and there is no recurrence.
So the number of the title of this paper `almost' made it to the OEIS, since it is the $n=13$ entry of that sequence.
The sequence corresponding to $k=5$ is {\tt A123297}, there is no recurrence and it only goes as far as $n=11$.
The sequences for $k \geq 6$ were not present. The above output file has many terms, and {\bf recurrences} for {\it all} these sequences
through $k=19$.

What about the sequences $\{ F[n](k)\}_{k=0}^{\infty}$ for a fixed, numeric $n$? Now things are much more complicated (and slower),
and the Almkvist-Zeilberger algorithm is not applicable. While we know, {\it a priori} that such a recurrence
{\bf exists}, for {\it every} $n$, once $n$ gets larger, it becomes computationally challenging. The Mathematica package [K], written by one of us (CK),
can handle it very well, and we found linear recurrences for $n \leq 9$. 

The output file
{\tt https://sites.math.rutgers.edu/\~{}zeilberg/tokhniot/oMultiDer2.txt} \quad 
contains such recurrences for $2 \leq n \leq 9$.

Note that the case $n=2$ is the identically $1$ sequence (after all, the Laguerre polynomials are orthonormal,
and there exists {\bf exactly one} multi-set derangement of $1^n2^n$). The case $n=3$ gives the
Franel sequence $\sum_{k=0}^{n} {{n} \choose {k}}^3$, as noted by Richard Askey (mentioned at the end of [EG],
and in [As], p.43). Surprisingly the sequences for $n \geq 4$ were not in the OEIS.

{\bf Conclusion}

People have been playing cards for centuries, and luminaries such as Persi Diaconis analyzed their random shuffling.
It is surprising that such a natural number, the number of ways of deranging a standard deck, ignoring suits, could not
be found anywhere in the OEIS, or for that matter in the internet. In fact, an internet search for the number of
the title did give one hit, but unless you know the number beforehand, it is useless.

But the main message of this article, in addition to filling this {\it much needed gap}, is to illustrate the great
utility of the Even-Gillis formula that, interfaced with symbolic-computation (the holonomic systems approach [AZ][Z][K])
can easily find , almost immediately, the number of ways of, say, deranging the multisets $1^{19} \dots 2000^{19}$ 
and the number of ways of deranging the multiset $1^{2000} \dots 9^{2000}$, both of which are contained in the above-mentioned
output files.

{\bf References}

[AlZ] Gert Almkvist and Doron Zeilberger, {\it The method of differentiating under the
integral sign}, J. Symbolic Computation {\bf 10}, 571-591 (1990). \hfill\break
{\tt https://sites.math.rutgers.edu/\~{}zeilberg/mamarim/mamarimhtml/duis.html} \quad .

[As] Richard Askey, {\it Orthogonal Polynomials and Special Functions}, SIAM, 1975.

[EG] S. Even and J. Gillis, {\it Derangements and Laguerre polynomials}, Math. Proc. Camb. Phil. Soc. {\bf 79} (1976), 135-143.

[K] Christoph Koutschan, {\it   Advanced applications of the holonomic systems approach}, 
PhD thesis, Research Institute for Symbolic Computation (RISC), Johannes Kepler University, Linz, Austria, 2009.\hfill\break
{\tt http://www.koutschan.de/publ/Koutschan09/thesisKoutschan.pdf}, \hfill\break
{\tt http://www.risc.jku.at/research/combinat/software/HolonomicFunctions/} \quad.

[M] Percy Alexander MacMahon, {\it Combinatory Analysis}, in: Encyclopedia Britannica, eleventh edition, 1910.
Volume VI, pp. 752-758.

[OEIS] {\it The On-Line Encyclopedia of Integer Sequences}, {\tt https://oeis.org} \quad .

[R] John Riordan, `Introduction to Combinatorial Analysis', Dover, originally published by John Wiley, 1958.

[Z] Doron Zeilberger, {\it A Holonomic systems approach to special functions
identities}, J. of Computational and Applied Math. {\bf 32}, 321-368 (1990). \hfill\break
{\tt  https://sites.math.rutgers.edu/\~{}zeilberg/mamarim/mamarimhtml/holonomic.html }\quad .

\bigskip
\hrule
\bigskip
Shalosh B. Ekhad, Department of Mathematics, Rutgers University (New Brunswick), Hill Center-Busch Campus, 110 Frelinghuysen
Rd., Piscataway, NJ 08854-8019, USA. \hfill\break
Email: {\tt ShaloshBEkhad at gmail  dot com}   \quad .
\bigskip
Christoph Koutschan, Johann Radon Institute of Computational and Applied Mathematics (RICAM), Austrian Academy of Sciences,
Altenberger Strasse 69, A-4040 Linz, Austria \hfill\break
Email: {\tt  christoph.koutschan at ricam dot oeaw dot ac dot at}   \quad .
\bigskip
Doron Zeilberger, Department of Mathematics, Rutgers University (New Brunswick), Hill Center-Busch Campus, 110 Frelinghuysen
Rd., Piscataway, NJ 08854-8019, USA. \hfill\break
Email: {\tt DoronZeil at gmail  dot com}   \quad .

\bigskip

{\bf Jan. 25, 2021} 

\end